\documentclass[11pt]{article}
\usepackage{amssymb, amsmath, amsthm, amsfonts,xcolor, enumerate}
\usepackage{fullpage}
\usepackage{hyperref}

\usepackage{t1enc}
\usepackage[T1]{fontenc}
\usepackage[utf8]{inputenc}

\newtheorem{theorem}{Theorem}[section]
\newtheorem{lemma}[theorem]{Lemma}

\newtheorem{defn}[theorem]{Definition}

\newtheorem{claim}[theorem]{Claim}
\newtheorem{prop}[theorem]{Proposition}
\newcommand{\ep}{\varepsilon}

\def\qedf{\hfill $\Box$}

\def\komment#1{}
\let\komment=\footnote

\begin{document}

\title{A new graph decomposition method\\ for bipartite graphs}

\author{
B\'ela Csaba}

\date{}

\maketitle
\centerline{\large Extended Abstract}

\begin{abstract}
Given a sufficiently large and sufficiently dense bipartite graph $G=(A, B; E),$ we present a novel method for decomposing the majority of the edges of $G$ into quasirandom graphs so that the vertex sets of these quasirandom graphs partition the majority of $A.$ The method works for relatively small or sparse graphs, and can be used to substitute the Regularity lemma of Szemer\'edi in some 
graph embedding problems.
\end{abstract}

\section{Introduction}
All graphs considered in this paper are simple. 
The Szemer\'edi Regularity lemma~\cite{Szemeredi} is one of the most powerful tools of graph\footnote{There are also
hypergraph versions that play crucial role in extremal hypergraph theory and combinatorial number theory, see eg.,~\cite{GowersHyp} or~\cite{RodlHyp}.} theory. 
It is also used in many areas outside graph theory, for example 
in number theory and algorithms. 

\begin{theorem}[Szemer\'edi]\label{RegLemma}
For every $\varepsilon>0$ there exists a $n_0=n_0(\varepsilon)>0$ such that if $G$ is a simple graph on $n\ge n_0$ 
vertices $G$ admits an $\varepsilon$-regular equipartition of its vertex set into at most $k=k(\ep)$ parts.
\end{theorem}

We will give a short introduction to the necessary notions in the next section. Here we only mention that $\varepsilon$-regularity is a notion of quasirandomness, and equipartition means, roughly, the partition of the vast majority of the vertex set of $G$ into equal sized subsets so that all, but an $\varepsilon$ proportion of the
pairs of subsets span an $\varepsilon$-regular bipartite subgraph of $G.$

The dependence of $n_0$ on $\varepsilon$ in Theorem~\ref{RegLemma} is determined by a {\it tower function} $T$ evaluated at $1/\varepsilon^5,$ where $T$ can be defined inductively as follows: $T(1)=2,$ and for $i>1$ we have $T(i)=2^{T(i-1)}.$ Hence, the value of $n_0$ makes the Regularity lemma essentially impractical.
It is also well-known that we cannot hope for a much better bound, since as was proven by Gowers~\cite{Gowers} and more recently by Conlon and Fox~\cite{CF}, there are graphs for which the number of clusters in the Regularity lemma is necessarily a tower function of $1/\varepsilon.$ Note also, that the lemma is only meaningful for so called dense graphs, that is, graphs that contain a constant proportion of the possible edges.

In this paper we present a new graph decomposition method for bipartite graphs, which can be applied for graphs of practical size and for graphs having vanishing density.
While the Regularity lemma is useful in many areas of mathematics and computer science, our contribution may not be so widely applicable. Still, it can be used for finding certain subgraphs in a host graph. As an illustration we will give the details of a tree embedding algorithm that uses this graph decomposition method. 
Let us mention that Gowers in~\cite{rank1} presented a decomposition for bipartite graphs that is somewhat similar to the one discussed here, and used it for a problem in number theory.
That decomposition has different parameters and a much longer and harder proof.

Due to the importance of the Regularity lemma several researchers found weakened versions (\cite{AlonDLRY},~\cite{FriezeKannan},  etc.) in which the dependence of $\varepsilon$ and $n_0$ is not determined by a tower function. These are important developments with several applications, still, none of them seems to be so widely applicable as the original one. 
The so called absorption method~\cite{Szemeredi2} is also a choice for avoiding the use of the Regularity lemma in some embedding problems. 

The outline of the paper is as follows. First, we provide the necessary notions for the decomposition and then describe the decomposition method in the next section. 
In the subsequent section we provide an application, namely, we show that we can find a large subtree in a graph on $n$ vertices having $\Omega(n^2 \log\log n/{\log n})$ edges. 

\section{Definitions, main result}\label{defek}

Given a graph $G$ with vertex set $V$ and edge set $E$ we let $deg_G(v)$ denote the degree of $v\in V.$ If it is clear from
the context,
the subscription may be omitted. The neighborhood of $v$ is denoted by $N(v),$ so $deg(v)=|N(v)|.$ The minimum degree of $G$ is denoted by $\delta(G).$ If $S\subset V,$ then $deg(v; S)=|N(v)\cap S|.$  The set of edges between two disjoint
sets $S, T\subset V$ is denoted by $E(S, T),$ and we let $e(S, T)=|E(S, T)|.$ We also let $e(G)=|E(G)|.$

Let $G=G(A, B; E)$ be a bipartite graph. The {\it density} $d_G(A, B)$ or, if $G$ is clear from the context, $d(A, B),$ is defined as follows:$$d(A, B)=\frac{e(G)}{|A|\cdot |B|}.$$   

Given a number $\varepsilon \in (0,1)$ we say that $G$ is an {\it $\varepsilon$-regular pair} if the following holds for every 
$A'\subset A,$ $|A'|\ge \varepsilon |A|$ and $B'\subset B,$ $|B'|\ge \varepsilon |B|$: $$|d_G(A, B)-d_G(A', B')|\le \varepsilon.$$ 
The $\varepsilon$-regular equipartition of a graph $G$ on $n$ vertices means that there exists a number $k=k(\ep)$ and $V(G)=V_0\cup V_1 \cup \ldots \cup V_k$ such that $V_i\cap V_j=\emptyset$ for $i\neq j,$ $|V_0|\le \varepsilon n,$ $||V_i|-|V_j||\le 1$ for every $1\le i, j \le k$ and all but at most $\varepsilon k^2$ 
pairs
$V_iV_j$ are $\varepsilon$-regular for $1\le i, j.$ The $V_i$ sets are called clusters, and $V_0$ is the exceptional cluster. 

Roughly speaking, the Regularity lemma asserts that every graph can be well approximated by a collection of quasirandom graphs that are defined between the non-exceptional clusters. Unfortunately, the number of non-exceptional clusters is a tower function of $1/\varepsilon.$

While our goal is to provide an alternative for the Regularity lemma, we will also make use of the regularity concept. Our definition is slightly more permissive than the usual one above, this enables us to
give a very short proof of our decomposition, and it is still powerful enough to be applicable in several embedding problems. It is called {\it lower regularity}, and is used by other researchers as well. 

\begin{defn}
Given a bipartite graph $G=G(A, B)$ we say that $G$ is a lower $(\varepsilon, \eta, \gamma)$-regular pair, if for any $A'\subset A, B'\subset B$ with $|A'|\ge \varepsilon |A|,$ $|B'|\ge \eta |B|$ we have $e(A', B')\ge \gamma \cdot |A'|\cdot|B'|.$
\end{defn}

Note that in the usual definition of an $\varepsilon$-regular pair one has $\varepsilon=\eta,$ and the edge density between two sufficiently large subsets is between $d_G-\varepsilon$ and
$d_G+\varepsilon.$ We want to have flexibility in this notion, and allow sub-pairs with relatively low density, and the $\varepsilon\neq \eta$ case, too.

We are ready to state our main result, the precise formulation is as follows.

\begin{theorem}\label{particios}
Let $G=G(A, B)$ be a bipartite graph  with vertex classes $A$ and $B$ such that $|A|=n$ and $|B|=m,$ and every
vertex of $A$ has at least $\delta m$ neighbors in $B.$ Let $0<\varepsilon, \eta, \gamma <1$ be numbers 
so that $\eta\le 1/6$ and $\gamma\le \min\{\eta/4, \delta/20\}.$ Then there exists a partition $A=A_0\cup A_1 \cup \ldots \cup A_k,$ and $k$ not necessarily disjoint subsets of $B,$ $B_1, \ldots, B_k$
such that $|A_i|\ge \varepsilon \cdot \exp\left(-2\log(\frac{1}{\varepsilon})\log(\frac{2}{\delta})/\eta\right)n$ 
for $i\ge 1,$ $|A_0|\le \epsilon n,$ the subgraphs $G[A_i, B_i]$ for $1\le i \le k$ are all lower $(\varepsilon, \eta, \gamma)$-regular, and
$$\sum_{i=1}^ke(G[A_i, B_i])\ge e(G)-(\varepsilon+2\gamma)nm.$$
Moreover, $$k\le \frac{2}{\varepsilon \delta} e^{2\log(\frac{1}{\varepsilon})\log(\frac{2}{\delta})/\eta}.$$
\end{theorem}


\section{Proof of Theorem~\ref{particios}}

Let us remark that we will not be concerned with floor signs, divisibility, and so on in the proof. This makes the notation simpler,
easier to follow.

As we have seen edge density plays an important role in regularity. We need a simple fact which is called {\it convexity of density} (see 
eg.~in~\cite{KS}), the proof is left for the reader.

\begin{claim}\label{convexity} Let $F=F(A, B)$ be a bipartite graph, and let $1\le k \le |A|$ and $1\le m \le |B|.$ Then 
$$d_F(A, B)=\frac{1}{\binom{|A|}{k}\binom{|B|}{m}}\sum_{X\in \binom{A}{k}, Y\in \binom{B}{m}} d(X, Y).$$
\end{claim}

In order to prove Theorem~\ref{particios} we need a lemma 
that is the basic building block of our decomposition method. 

\begin{lemma}\label{egylepes}
Let $F=F(A, B)$ be a bipartite graph  with vertex classes $A$ and $B$ such that $|A|=a$ and $|B|=b,$ and every
vertex of $A$ has at least $\delta b$ neighbors in $B.$ Let $0<\varepsilon, \eta, \gamma <1$ be numbers 
so that $\eta\le 1/6$ and $\gamma\le \min\{\eta/4, \delta/20\}.$ 
Then $F$ contains a lower $(\varepsilon, \eta, \gamma)$-regular pair $F[X, Y]$ such that $|X|\ge  \exp\left(-2\log(\frac{2}{\varepsilon})\log(\frac{2}{\delta})/\eta\right)a$ and $|Y|\ge (\delta(1 - \eta)-2\gamma)b.$ 
\end{lemma}


\noindent {\bf Proof:} We prove the lemma by finding two sequences of sets $X_0=A, X_1, \ldots, X_l$ and $Y_0=B, Y_1,$ 
\ldots, $Y_l$ such that  for every $1\le i\le l$ we have that $X_i\subset X_{i-1}$ and $Y_i\subset Y_{i-1}$ and $$\varepsilon |X_{i-1}|/2 \le |X_i|\le \varepsilon |X_{i-1}|$$ and $$|Y_i|=(1-\eta)|Y_{i-1}|,$$ moreover, the last pair  $F[X_l, Y_l]$ is lower $(\varepsilon, \eta, \gamma)$-regular. 
Hence, we may choose $X=X_l$ and $Y=Y_l.$

We find the set sequences $\{X_i\}_{i\ge 1}$ and $\{Y_i\}_{i\ge 1}$ by the help of an iterative procedure. This procedure stops in the 
$l$th step, if $F[X_l, Y_l]$ is lower $(\varepsilon, \eta, \gamma)$-regular. We have another stopping rule: if 
$|Y_l| \le (\delta (1+\eta/2)-2\gamma)b$ for some $l,$ we stop. Later we will see
that in this case we have found what is desired, $F[X_l, Y_l]$ must be a lower $(\varepsilon, \eta, \gamma)$-regular pair.

In the beginning we check if $F[X_0, Y_0]$ is a lower $(\varepsilon, \eta, \gamma)$-regular pair. If it is, we stop. If not then $X_0$ has a subset $X'_1$ precisely of size $\varepsilon |X_0|$ and $Y_0$ has a subset 
$Y'_1$ precisely of size $\eta|Y_0|$ such that $e(F[X'_1, Y'_1])<\gamma |X'_1|\cdot |Y'_1|,$ here we used Claim~\ref{convexity} in order to obtain the sizes of $X'_1$ and $Y'_1.$

 Let $X''_1$ be the set of those vertices of $X'_1$ that have
more than $2\gamma |Y'_1|$ neighbors in $|Y'_1|.$
Simple counting shows that $|X''_1|\le |X'_1|/2.$ Let $X_1=X'_1-X''_1,$ those vertices of $X'_1$ that have less than 
$2\gamma |Y'_1|$ neighbors in $|Y'_1|.$ By the above we have $|X'_1|/2\le |X_1|\le |X'_1|.$ Set $Y_1=Y_0-Y'_1.$  

For $i\ge 2$ the above is generalized. If $F[X_{i-1}, Y_{i-1}]$ is not  a lower $(\varepsilon, \eta, \gamma)$-regular pair then we do the following. First find $X'_i\subset X_{i-1}$ 
and $Y'_i\subset Y_{i-1}$ such that $|X'_i| =\varepsilon |X_{i-1}|$ and $|Y'_i| =\eta |Y_{i-1}|$ and $e(F[X'_i, Y'_i])<\gamma |X'_i|\cdot |Y'_i|.$ Similarly to the above we define $X_i\subset X'_i$ to be the set of those vertices of $X'_i$ that have less than $2\gamma |Y'_i|$ neighbors in $Y'_i.$ As before, we have $|X'_i|/2\le |X_i|\le |X'_i|.$ Finally, we let $Y_i=Y_{i-1}-Y'_i.$

Using induction one can easily verify that the claimed bounds for $|X_i|$ and $|Y_i|$ hold for every $i.$ It might not be so clear that this process stops in a relatively few iteration steps.  
For that we first find an upper bound for the number of edges that connect the vertices of $X_i$ with $B-Y_i.$ If $u\in X_i$ then $u$ have at most $2\gamma (|Y'_1|+\ldots +|Y'_i|)\le 
2\gamma b$ neighbors in $B-Y_i$ using that $Y'_s\cap Y'_t =\emptyset$ for every $s\neq t.$

Next we show that if $(\delta(1+\eta/2)-2\gamma)(1-\eta)b< |Y_l|\le (\delta(1+\eta/2)-2\gamma)b$ then $F[X_l, Y_l]$ must be lower regular. 
Assume that $u\in X_l.$ Then $deg(u; Y_l)\ge (\delta-2\gamma)b,$ using our argument above, hence, the number of 
{\it non-neighbors} of $u$ in $Y_l$ is at most $(\delta(1+\eta/2)-2\gamma)b-(\delta -2\gamma)b=\delta\eta b/2.$ 
Let $Y'\subset Y_l$ be arbitrary with $|Y'|=\eta |Y_l|.$ Then 
$$\frac{5}{6}\eta(\delta(1+\eta/2)-2\gamma)b\le |Y'|\le \eta(\delta(1+\eta/2)-2\gamma)b,$$ using that $\eta\le 1/6.$
We have
$$deg(u; Y')\ge |Y'|-\delta\eta b/2\ge \frac{5}{6}\eta( \delta(1+\eta/2)-2\gamma)b -\delta\eta b/2.$$
Using the upper bounds we imposed on $\eta$ and $\gamma,$ one easily obtains that 
$$deg(u; Y')\ge (\delta\eta/3+5\delta\eta^2/12-5/3\gamma \eta)b\ge \gamma |Y'|.$$

Hence, for every $X'\subset X_l$ and $Y'\subset Y_l$ with $|Y'|=\eta |Y_l|$ we have $$e(X', Y')\ge \gamma |X'|\cdot |Y'|,$$  
that is, if the procedure stopped because we applied the stopping rule, then the resulting pair must always be lower $(\varepsilon, \eta, \gamma)$-regular. Of course, this means that no matter how the procedure stops, it finds a lower regular pair.

Next we upper bound the number of iteration steps. In every step the $Y$-side shrinks by a factor of $(1-\eta).$ 
We also have that $|Y_l|> (\delta(1+\eta/2)-2\gamma)(1-\eta)b.$ Putting these together we get that
$$(1-\eta)^l>(\delta(1+\eta/2)-2\gamma)(1-\eta)>\delta/2.$$ Hence, $$l< \frac{\log (2/\delta)}{\log(1/(1-\eta))}<2\frac{\log(2/\delta)}{\eta},$$
here we used elemantary calculus (in particular, the Taylor series expansion of $\log(1+x)$) and our condition that $\eta$ 
is less than 1/6.

What is left is to show the lower bound for $|X_l|.$ Note, that $|X_i|/|X_{i-1}|\ge \varepsilon/2$ for every $i\ge 1.$ Hence, 
$$|X_l| \ge \left(\frac{\varepsilon}{2}\right)^la =e^{-2\log(2/\varepsilon)\log(2/\delta)/\eta}a.$$
 
\qedf

\medskip

We are ready to prove the main result of the paper.

\smallskip

\noindent {\bf Proof} (of Theorem~\ref{particios}):
The proof is based on iteratively applying Lemma~\ref{egylepes}. First we apply Lemma~\ref{egylepes} for $G$ and find
a lower $(\varepsilon, \eta, \gamma)$-regular pair $G[X_l, Y_l],$ where $X_l\subset A$ and $Y_l\subset B.$ Let $A_1=X_l$
and $B_1=Y_l.$ 
Next we repeat this procedure for the graph $G[A-A_1, B].$ Similarly to the above we define the $A_2$ and $B_2$ sets, where
$A_2\subset A-A_1,$ $B_2\subset B,$ and $G[A_2, B_2]$ is a lower $(\varepsilon, \eta, \gamma)$-regular pair. 

Continue this way, finding the lower regular pairs $G[A_i, B_i]$ using Lemma~\ref{particios} such that 
$A_i\subset A-(A_1\cup \ldots \cup A_{i-1}),$ $B_i\subset B,$ and $G[A_i, B_i]$ is a lower $(\varepsilon, \eta, \gamma)$-regular pair. We stop when $$|A-(A_1\cup \ldots \cup A_{i})|<\varepsilon |A|.$$ At this point set 
$A_0 = A-(A_1\cup \ldots \cup A_{i}).$ 

Let us now prove the upper bound for the number of pairs in the decomposition. As we have shown earlier $|A_i|\ge \exp\left(-2\log(\frac{2}{\varepsilon})\log(\frac{2}{\delta})/\eta\right)n$ for $i\ge 1.$ The number of edges in an $A_iB_i$
pair is at least $|A_i|(\delta -2\gamma)m>|A_i|\delta m/2.$ For any $1\le i\neq j\le k$ the edge sets of the pairs
$A_iB_i$ and $A_jB_j$ are disjoint, and the total number of edges in lower regular pairs is at most $nm.$ Hence, we have
$$k\le \frac{2nm e^{2\log(\frac{1}{\varepsilon})\log(\frac{2}{\delta})/\eta}}{\varepsilon \delta nm}=
\frac{2}{\varepsilon \delta} e^{2\log(\frac{1}{\varepsilon})\log(\frac{2}{\delta})/\eta}.$$

There is only one question left, bounding the total number of edges that belong to the lower regular pairs. Assume first that
$u\in A-A_0.$ We saw earlier in Lemma~\ref{particios} that $u$ lost at most $2\gamma|B|$ edges. This explains 
the $2\gamma mn$ term in the theorem. If $u\in A_0,$ none of the edges incident to it belongs to any of the
lower regular pairs, however, $|A_0|\le \varepsilon n,$ therefore, the total number of edges incident to vertices of $A_0$
is at most $\varepsilon nm.$ With this
we found the decomposition of $G$ what was desired.
\qedf


Let us finish this section with a remark. Without the lower bound for the sizes of the $A_i$ sets, Theorem~\ref{particios} would
be trivial: every vertex $v\in A$ could be a ``subset'' $A_v$ (a singleton), and its neighborhood $N(v)$ is the corresponding $B_v.$
The result is interesting only when the $A_i$ sets are large. For example, let $G$ be the following. It is a sparse bipartite graph   
with vertex classes $A$ and $B$ such that $|A|=|B|=n.$ Set $\varepsilon = \eta= 1/10,$ $\delta=\log\log n/\log n,$ and $\gamma=\delta/20.$ Then $G$ has $O(n^2\log\log n/\log n)$ edges, and the $A_i$ sets for $i\ge 1$ have size
$\Omega(n/(\log n)^c),$ where $c< 60,$ and every $(A_i, B_i)$ pair is a lower $(0.1, 0.1, \log \log n/(20 \log n))$-regular pair. 

\section{An application}

The main advantage of Theorem~\ref{particios} is that, as the above example shows, it can be applied for graphs 
having ``real-life'' size, or for relatively sparse graphs, unlike the Szemer\'edi Regularity lemma.  
Therefore, it may extend the scope when usual methods for graph embedding (eg. counting lemma or the Blow-up lemma~\cite{KSSz}) can be applied. 


Given a tree $T$ rooted at $r$ its {\it level sets} are defined as follows: $L_1=r,$ $L_2=N(r),$ in general, $L_{i+1}=N^i(r),$ etc., where 
$N^i(r)$ denotes those vertices of $T$ that are exactly at distance $i$ from $r$ in $T.$

\begin{prop}\label{fa}
Let $0< \varepsilon, \eta, \gamma <1/10$ such that $\eta=4\gamma$ and $\varepsilon=\gamma^2/10.$ 
Assume $G[A, B]$ is a lower $(\varepsilon, \eta, \gamma)$-regular pair. Let $T$ be a tree rooted at $r,$ having color 
classes $X$ and $Y$ such that $r\in X,$ $|X|\le (1-10\gamma)|A|$ and $|Y|\le (1-10\gamma)|B|.$ Assume further that 
for every $i\ge 1$ we have $|L_{2i}|\le \varepsilon |A|$ and $|L_{2i+1}|\le \eta |B|.$
Then $T\subset G[A, B].$  
\end{prop} 

Let us remark that $T$ does not have to have bounded degree, unlike in many tree embedding results. In fact, it can have vertices with linearly large degrees, if $\delta$ and the other parameters are constants. 
The statement holds for every $G$ for which Lemma~\ref{egylepes} can be applied, hence, $G$ can have $o(n^2)$ edges.

We need the following simple claim, the proof is left for the reader.

\begin{claim}\label{kovetkezmeny}
Let $F=F(U, V)$ be a lower $(\varepsilon, \eta, \gamma)$-regular pair. Let $U'\subset U$ and $V'\subset V$ such that 
$|U'|\ge \varepsilon |U|$ and $|V'|\ge \eta |V|.$ Then $U'$ can have at most $\varepsilon |U|$ vertices that have less
than $\gamma |V'|$ neighbors in $V'.$ Similarly, $V'$ can have at most $\eta |V|$ vertices that have less than $\gamma |U'|$
neighbors in $U'.$  
\end{claim}


\noindent {\bf Proof of the proposition:} We prove via an embedding algorithm. Let $X=\{x_1, \ldots, x_k\}$ and 
$Y=\{y_1, \ldots, y_m\},$ where $r=x_1.$ We will find the images of the vertices of $T$ so that we embed height-2 subtrees of $T$ in every step, having vertices from $Y$ in the middle level. 

Denote $\varphi: V(T)\longrightarrow A\cup  B$ the edge-preserving mapping that we construct.
Let $A^f,$ respectively, $B^f$ denote the {\it free} (ie. vacant) vertices of $A,$ respectively, $B.$ 
These sets are shrinking as the embedding of $T$ proceeds, but due to
the conditions of Proposition~\ref{fa} we always have that $|A^f|\ge 10\gamma |A|$ and $|B^f|\ge 10\gamma |B|.$
Divide $A^f$ randomly into three disjoint, approximately equal-sized subsets $A^f_1, A^f_2$ and $A^f_3.$ Let $B_1'\subset B^f$
be the set of those vertices that have less than $\gamma |A|$ neighbors in $A^f_1,$ the sets $B_2'$ and $B_3'$ are defined
analogously. Then $|B_1'|, |B_2'|, |B_3'|\le \eta |B|.$

Let $v$ be an arbitrary vertex of, say, $A^f_1$ that has at least $\gamma |B^f-B_1'-B_2'-B_3'|$ neighbors in $B^f-B_1'-B_2'-B_3'.$ By Claim~\ref{kovetkezmeny} we know that $A^f_1$ has many such vertices. By the definition of the $B_i'$ sets we have that
every vertex in $N(v)$ has at least $\gamma |A^f_i|/4$ neighbors in $A^f_i$ for $i=1, 2, 3.$ Pick the largest of the $A^f_i$
sets, say, it is $A^f_2.$ Then the height-2 subtree originating at $r$ will be embedded so that the neighbors of $r$ will be mapped onto $N(v)$ arbitrarily ($|L_2|$ is smaller, than $|N(v)|$), and by construction every vertex of $N(v)$ will have many neighbors
in $A^f_2.$ Now we redetermine the subsets $B_1', B_2', B_3',$ as some vertices have become covered in $A$ and in $B.$ 
For the third level of the height-2 subtree originating at $r$ we take those vertices of $A^f_2$ that are neighboring with
at least a $\gamma$ proportion of $B^f-B_1'-B_2'-B_3'.$  
Note that for every $\phi(y)$ where $y$ is in the middle level we have many choices: except at most $\varepsilon |A|$ vertices 
of $A^f_2$ the neighborhood $N(\phi(y))$ contains vertices with large degrees into $B^f-B_1'-B_2'-B_3'.$ 
This means that we are able to map the third level. 
Next we continue this process so that we embed the height-2 subtrees originating at the vertices of the third level one-by-one. 

There is only one missing detail here, the reason why we divided $A^f$ randomly in the beginning: if we have three $A^f_i$ sets, then the {\it active level} belongs to one of them, say, it is $A^f_i.$ Then we map the vertices of $T$ that are exactly two levels below them into the larger $A^f_j$-set, where $j\in \{1, 2, 3\}-i.$ This way we never eat up any of the $A^f_i$ sets at any
point in time. Since the color classes of $T$ are sufficiently small, this procedure never gets stuck.
\qedf

\medskip

We have just showed how to embed an almost spanning tree into {\it one} lower regular pair. This can be used
to approximately tile the edge set of a sufficiently dense (say, having $\Omega(n^2\log\log n/\log n)$ edges and $n$ vertices) graph $G$ by large edge-disjoint subtrees. 
The rough sketch of this approximate decomposition is as follows. Apply Theorem~\ref{particios} for the graph $G,$ and then
using Proposition~\ref{fa} find one-one almost spanning subtree in the lower regular pairs. Delete the edges used for the
subtrees. If the resulting graph has sufficiently many edges then one can use Theorem~\ref{particios} again, and then
Proposition~\ref{fa} for every lower regular pair. The process stops when the remaining vacant subgraph of $G$ is too
sparse, and therefore one cannot find many large degree vertices in it. Hence, with this method one can tile the vast majority
of edges of a graph having sufficiently large density by large subtrees. We leave the details for the reader.




\end{document}